\documentclass{tip}
\usepackage{orcidlink}
\newcommand{\autora}{R.~Zöllner}
\newcommand{\autorb}{F.~Schuricht}
\newcommand{\autorc}{T.~Schmidt}
\newcommand{\autord}{W.~Hofmann}
\title{\vspace{-2.5cm}\hrulefill\\\color{darkblue}\sffamily\LARGE\bfseries Semi-analytical Approach to Trajectory Optimization for Stacker Cranes Regarding Energy Saving\\}
\author{\sffamily\bfseries \autora\orcidlink{0000-0002-3544-6622}${}^a$,~\autorb${}^b$,~\autorc\orcidlink{0000-0002-1484-7187}${}^a$~\&~\autord${}^c$}
\affil{\sffamily${}^a$Institute of Material Handling and Industrial Engineering, TU~Dresden,\newline01062 Dresden, Germany}
\affil{\sffamily${}^b$Faculty of Mathematics, TU~Dresden,
	01062 Dresden, Germany}
\affil{\sffamily${}^c$ Institute of Electrical Power Engineering, TU Dresden,	\newline01062 Dresden, Germany}
\date{\sffamily\today\\\hrulefill}

\begin{document}
	\maketitle
	\thispagestyle{plain}
	\vspace{-2.0em}
	
	\parbox[t]{0.30\linewidth}{\textbf{Keywords}
		\begin{flushleft}
			\begin{itemize}
			\item Trajectory Optimization
			\item Energy Recuperation
			\item Nonsmooth Analysis
			\item Warehouse Logistics
			\item Stacker Cranes 
				\end{itemize}
			\end{flushleft}}
	\hfill
	\parbox[t]{0.60\linewidth}{\textbf{Abstract}\\{\footnotesize The aim of this study is to give insights into the trajectory optimization w.r.t. energy consumption and recuperation for stacker cranes in a high-bay warehouse. Based on an analytical necessary optimality condition, a targeted numerical implementation is set up to perform systematic computations of optimal trajectories which are further categorized. Particularly, the differences between energy consumption and recuperation as well as for up and down movements are pointed out.
			
		Although examined for a concrete, experimentally validated model of stacker cranes, the methodical approach could be adapted to other electrical machines possesing a power flow model, i.e. a functional relation between the kinematics (velocity, acceleration for instance) and the resultant power. In addition, boundaries of the velocity, the acceleration and the jerk are incorporated.
		
		Such a systematic analysis of energy optimal trajectories can be further used for improving the job scheduling in a warehouse.
	
} }
	
	\vspace{2.0em} 
	\noindent\rule{\textwidth}{0.4pt}
	\vspace{-3.0em}
	
	\section{Introduction}\label{Sec1}
	Although energy optimization has been an active research field for several years, not all of its potential has yet been fully exploited. This case study focuses on saving electrical energy in the context of storage and retrieval machines (stacker crane) in high-bay warehouses. Due to their size and mass, the operation of stacker crane represent a key lever for a sustainable and efficient warehouse management. Typically, a stacker crane has two electrical drives: a running gear and a lifting gear, which can recuperate energy, for example, via an intermediate circuit. This allows energy to be saved, on the one hand, by drawing as little energy as possible from the grid, and, on the other hand, by optimally using the available energy for recuperation.
	
	Our paper addresses both objectives within the context of optimal trajectory planning for stacker cranes. As a first step, a model is developed from which the respective optimality condition is derived using variational principles (Section~\ref{Section2}). This allows the targeted implementation of a solution method (Section~\ref{NumStrat}), which, thanks to the previous analytical considerations, works significantly faster than would be the case with a direct, fully discretized method. A major advantage is the use of a non-equidistant time grid with few grid points (instead of an a priori discretized equidistant time grid). The power consumption and the power output of the two power trains of a stacker crane as the basis for energy considerations is discussed in Section~\ref{TechnExpAsp}. Section~\ref{Section3} presents the results of the numerical study.
	
	A comprehensive review of the immense existing literature would have been beyond the scope of this paper. For relevant review articles, please refer to~\cite{rouwenhorst,gu} systematizing warehouse design, operation and control, or refer to~\cite{Bartolini} reviewing green warehousing strategies. The papers~\cite{Jerman,Rucker,Faccio,Wu1} specifically address trajectory optimization and energy efficiency for stacker cranes (where~\cite{Jerman} considers energy regeneration, and~\cite{Rucker} deals with efficiency determination) and cranes (\cite{Wu1} focuses on optimal motion planning) -- all in all closely related, but slightly different approaches. Indirect and analytical solution methods are discussed in~\cite{Passenberg} (considering trucks and fuel consumption), and~\cite{Rams,Schwarzkopf,Zhang}, among others. Particularly, a variational approach to optimal motion planning for stacker cranes is presented in~\cite{Rams}, which -- however -- is more aligned to mechanical energies and oscillations of the stacker crane. Comparable approaches with technical considerations are provided, for instance, in~\cite{Pan,Li}.
		\section{Model Setup} \label{Section2}
\begin{figure}[h!]
		\centering
		\includegraphics[width=0.5\textwidth]{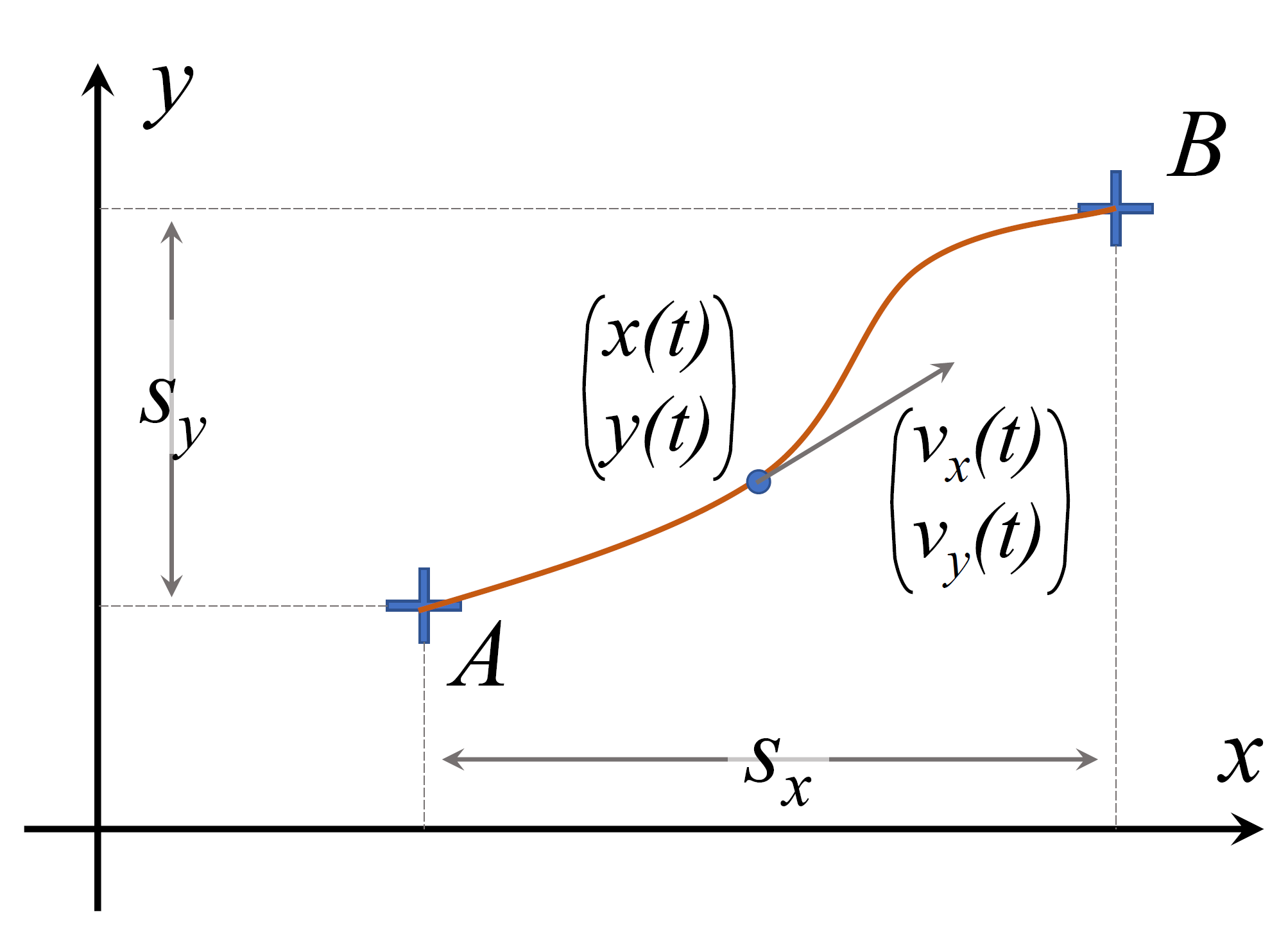}
		\caption{Geometric aspects of the trajectory optimization, where the load attachment device starts from point $A(x_A,y_A)$ and moves to point $B(x_B,y_B)$, in such a way that one of the drives -- exactly that one which needs more time to overcome the prescribed distance -- behaves time minimally and the other one makes optimal use of the energy recuperation~(\ref{3.1}) or minimizes the energy consumption~(\ref{3.11}). The respective distance of this drive is assigned to the parameter $s_0$ in~(\ref{3.2}).} \label{Figure1}
	\end{figure}

	The objective is either to maximize the energy recuperation between the two electrical drives of a stacker crane (see objective function~(\ref{3.1})) or to minimize the energy consumption~(see objective function~(\ref{3.11})) for a prescibed movement from point $A$ to a point $B$ in a vertical plane (see Fig.~\ref{Figure1}). Due to technical reasons, the velocity $v$, the acceleration $a=\dot{v}$ and the jerk $j=\dot{a}=\ddot{v}$ are bounded. In addition, the rectangle cornered by $A$ and $B$ must not be left and the functions $x(t)$ and $y(t)$ should be monotone. Due to logistical reasons, the energy saving strategy should not reduce the throughput. Therefore, the time minimal movements for both drives are evaluated. Such movements achieve the maximal velocity, acceleration and jerk (if the distance is long enough, otherwise only the acceleration and the jerk or only the jerk reaches their maximal values). This leads to a time duration $T_x$ for the horizontal movement ($T_y$ analogously). The available time horizon of the entire movement is defined as $T=\max \left\lbrace T_x, T_y\right\rbrace$, and consequently, one drive behaves time minimally (related quantities are indexed with ``slow''), and the other one (related quantities without index) is adjusted according to the following optimization problem.
	Let $P_{\mathrm{slow}}(t)$ denote the power profile of the drive which moves time minimally. And let $P(v,\dot{v})$ be a given sufficiently smooth function which models the power flow of the other~drive.

\newpage
\noindent	\textit{Objective: Energy Recuperation}
	
	The usage of recuperated energy is measured by
	\begin{equation}
		\min E_{\mathrm{rec}}(v):= \int_{0}^{T} |P_{\mathrm{slow}}(t) + P (v,\dot{v})|\, \mathrm{d}t.  \label{3.1}
	\end{equation}
	This objective function counts the amount of energy which is not recuperated where the absolute value penalizes energy flow in both directions (energy demand and recovery).
	\\
	
\noindent	\textit{Objective: Energy Consumption}
	
	The net energy consumption is measured by 
	\begin{equation}
		\min E_{\mathrm{con}}(v):= \int_{0}^{T} P_{\mathrm{slow}}(t) + P (v,\dot{v})\, \mathrm{d}t  \label{3.11}
	\end{equation}
	where energy demand and recovery are offset against each other.
	\\
	
\noindent	\textit{Constraints and Boundary Conditions}
	
	Independent of the objective function (\ref{3.1}) or (\ref{3.11}), one has the following constraints:
	\begin{align}
		G(v):&= \int_{0}^{T} v(t)~\mathrm{d}t-s_0=0 \label{3.2}, \\
		g_v(v):&= \max\limits_{t \in [0,T]} \left|v(t) \pm \frac{1}{2} v_{\mathrm{max}} \right| -\frac{1}{2}v_{\mathrm{max}} \leq 0  \label{3.3},\\ 
		g_a(v):&= \max\limits_{t \in [0,T]} | \dot{v}(t) | - a_{\max} \leq 0 \label{3.4},\\
		g_j(v):&= \max\limits_{t \in [0,T]} | \ddot{v}(t) | - j_{\max} \leq 0,\label{3.5}
	\end{align}
	supplemented by the boundary conditions 
	\begin{equation}
		v(0)=v(T)=\dot{v}(0)=\dot{v}(T)=0.
	\end{equation}
	Note that
	\begin{itemize}
		\item The parameter $s_0$ is either equal to $s_x$ or $s_y$, depending on which direction $x$ or $y$ belongs to the drive whose velocity profile is to be optimized (running gear or lifting gear).
		\item Condition (\ref{3.3}) ensures also the monotonicity of $x(t)$ or $y(t)$. 
		\item The combination of (\ref{3.2}, \ref{3.3}) impedes that the stacker crane leaves the admissible rectangle.
		\item The $\pm$ sign in (\ref{3.3}) indicates whether it is an up/down or left/right travel.
	\end{itemize}
	In more detail, if the velocity profile of the running gear is to be optimized (i.e. $s_0=s_x$) and $x_B > x_A$ (i.e. right travel) then~(\ref{3.3}) holds with a plus sign. This is also the case, if the lifting gear is subject to the optimization (i.e. $s_0=s_y$) and $y_B > y_A$ (up travel). Conversely, the minus sign applies if $s_0=s_x$ and $x_A > x_B$ (left travels) or $s_0=s_y$ and $y_A > y_B$ (down travels).
	\\
	
\noindent	\textit{Optimality Condition}
	
	According to nonsmooth variational first principles, there are two possible types of dynamics for $t\in\left] t_1,t_2\right[ \subseteq[0,T]$:
	\begin{itemize}[left=0pt, leftmargin=1cm]
		\item[CD)] Constraint Dominated Dynamics: There are subintervals, where exactely one of the constraints (\ref{3.3}~-~\ref{3.5}) is active, i.e.
		\begin{align}
			v(t)= \begin{cases}
				\pm v_\mathrm{max}\, \mbox{or}\, 0 &\mathrm{if \,(\ref{3.3}) \,is\, active}\\
				\pm a_\mathrm{max}t + v_0 &\mathrm{if \,(\ref{3.4}) \,is\, active}\\
				\pm \frac{1}{2}j_\mathrm{max}t^2 + a_0t + v_0 &\mathrm{if\, (\ref{3.5})\, is \,active}
			\end{cases}
		\end{align} 
		\item[EL)] Euler-Lagrange Dynamics: If no constraint is active the respective part of the trajectory is a solution of the Euler-Lagrange equation
		\begin{equation}
			0 = \frac{\partial P}{\partial v} - \frac{\mathrm{d}}{\mathrm{d}t} \frac{\partial P}{\partial \dot{v}} + \lambda_{G} \label{26}
		\end{equation}
		with a Lagrange multiplier arising from Equation~(\ref{3.2}).
	\end{itemize}
	If the objective function is (\ref{3.1}) there could theoretically occur a third type, namely 
	\begin{itemize}[left=0pt, leftmargin=1cm]
		\item[PB)] Power Balance Dynamics: The solution is given by the ODE
		\begin{equation}
			P_{\mathrm{slow}}(t) + 	P(v(t), \dot{v}(t)) = 0.
		\end{equation}
	\end{itemize}
	which is too restrictive and therefore irrelevant.
	
	In general, the task is now to find a time grid $(t_i)_{i=0}^{n}$ with $0 = t_0 < t_1 < ...< t_n = T$ and the respective type of dynamics for each subinterval $\left] t_i, t_{i+1}\right[ $ where all parameters as well as the grid points have to be chosen such that the entire function $v(t)$ is piecewise of class $\mathcal{C}^2$ and the boundary conditions are fulfilled. Thus, the problem remains as a low dimensional nonlinear optimization problem.
	\section{Numerical Strategy} {\label{NumStrat}}
	Beyond the above analysis, the remaining task can be considered as a continuous, convex and rather low dimensional optimization task. For such problems Matlab, Gurobi or Python provide solver, where we used scipy.optimize. The CD and EL dynamics have at most two free parameters (per time interval), and for a time grid $(t_i)^n_{i=0}$ of size $n+1$ we arrive at maximal $3n$ dimensions.
	
	Besides some auxiliary computations (see below), the main part of the implementation consists of an overarching iteration for the number of grid points $n+1$, where for $n$ fixed the respective optimum is calculated and $n$ is successively increased until the optimum could not be further improved. For fixed $n$, there is a two-step optimization. Within the first step, the function $P$ is approximated by a quadratic polynomial such that (\ref{26}) has an analytic solution and $v(t)$ and $E_{\mathrm{con}}$ as well as $E_{\mathrm{rec}}$ can be expressed with closed formulas. The solution of the first step serves as initial solution for the second step, where the approximation of $P$ is replaced by $P$ itself. The solution of (\ref{26}) as well as the interpretations for (\ref{3.1}) or (\ref{3.11}) are done numerically (using scipy.integrate).
	
\begin{figure}[!h]
		\centering
		\includegraphics[width=0.55\textwidth]{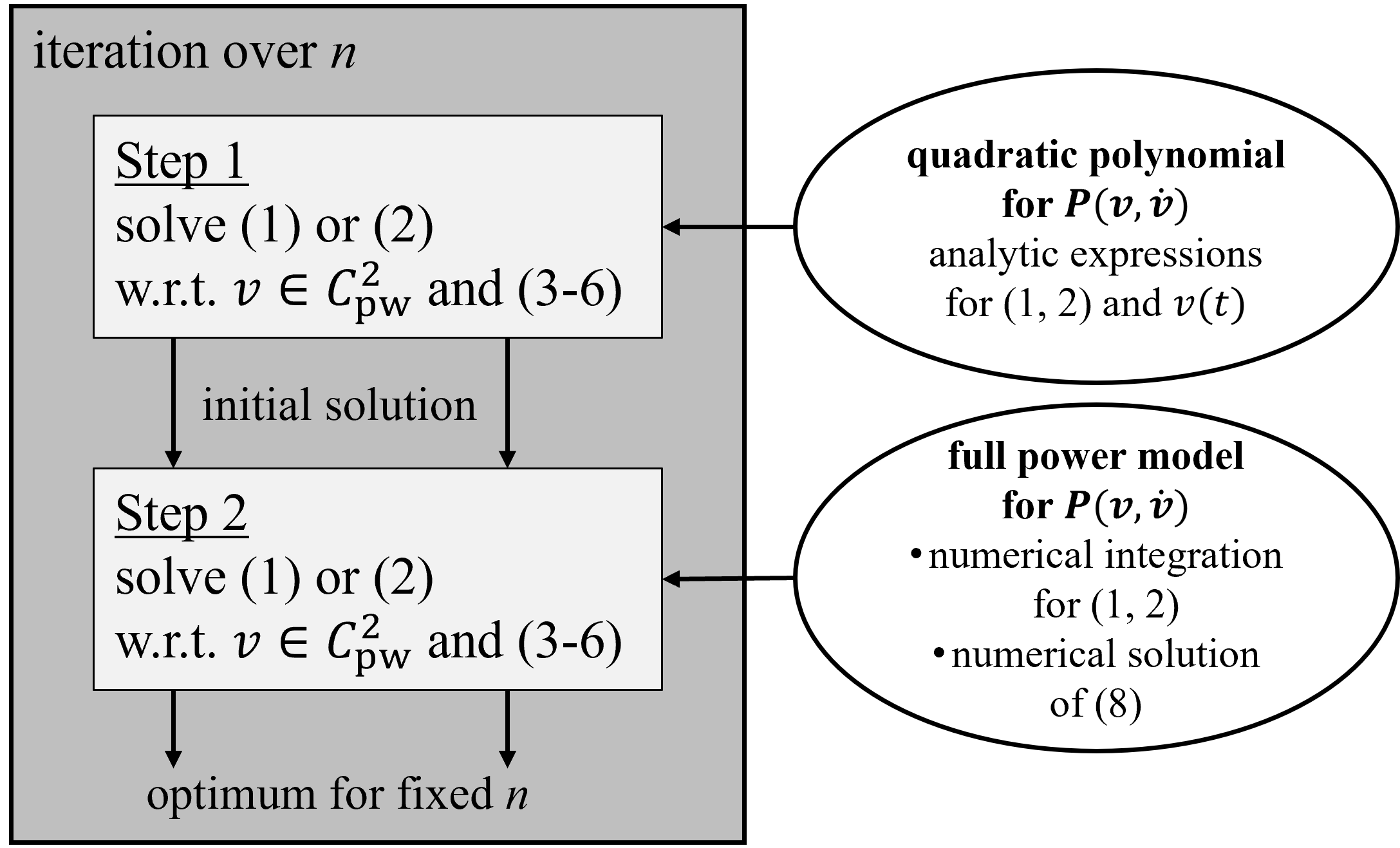} 
		\caption{Scheme of the implemented optimization strategy.}\label{iteration_over}
	\end{figure}

	Note that
	\begin{itemize}
		\item Since $P$ is assumed to be smooth, the solution of (\ref{26}) depends continuously (or even smoother) on the initial conditions and parameters.
		\item In addition, the initially gained solution of step 1 does not differ so much from the exact solution of step 2 since the solution of (\ref{26}) for a quadratic polynomial and the full $P$ does not differ very much.
		\item The fact that we are interested in piecewise $\mathcal{C}^2$ solutions leads to additional side conditions ensuring continuity (such as $v(t^-_i) = v(t^+_i)$ and $\dot{v}(t^-_i) = \dot{v}(t^+_i)$ for all inner grid points $i=1\ldots n-1$), which significantly reduces the search space. 
	\end{itemize}
	Figure~\ref{iteration_over} displays a scheme of the implementation.
	\section[Techn. and Exp. Aspects]{Technical and Experimental Aspects}\label{TechnExpAsp}
\begin{figure}[htbp!]
		\centering
		\subfloat[]{\includegraphics[width=0.45\textwidth]{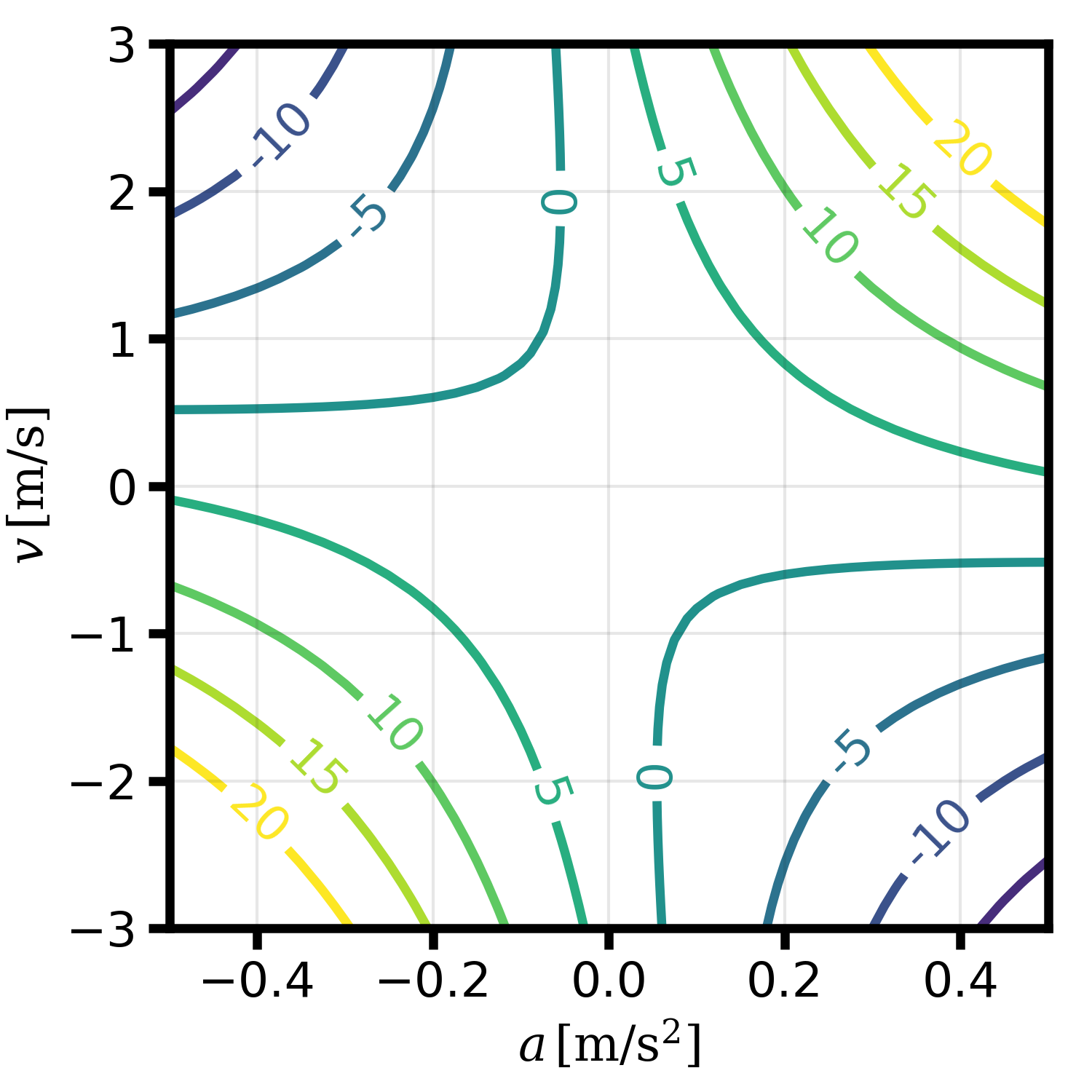}}
		\hfil
		\centering
		\subfloat[]{\includegraphics[width=0.45\textwidth]{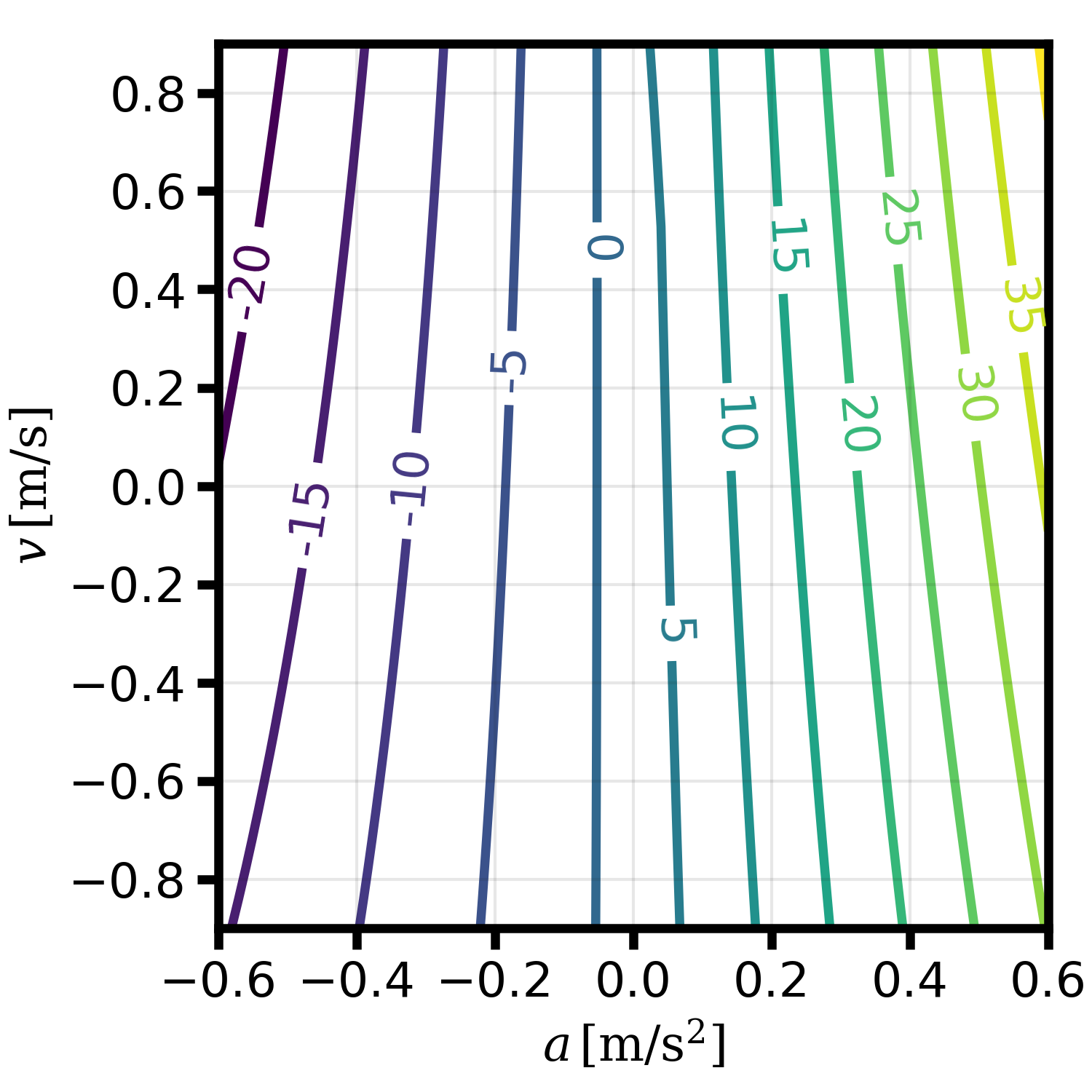}}
		\caption{Contour plots of the power $P=P(v,a)$ in kW as a function of the velocity and the acceleration for a load mass 1000~kg (panel (a): running gear, panel (b): lifting gear). Hereby, negative power means power supply.}
		\label{FigureB1}
	\end{figure}

	For any optimization of energy it is crucial to know the relation between the kinematic status of the vehicle (described by $v$ and $\dot{v}$) and resulting power $P$, i.e. the function $P(v,\dot{v})$. Such so called power flow models are discussed e.g. in \cite{Schutzhold,Cardenas,Lerher,Meneghetti,Monti,MM}, where we employed the model of \cite{Schutzhold}. Besides the pure mechanics, this model includes the gear, the asynchronous motor and the inverter as well as losses due to winding, iron, friction, as well as conduction and switching losses. We refrain here from a detailed technical description, instead Figure~\ref{FigureB1} shows the power as a function of the velocity~$v$ and the acceleration~$a=\dot{v}$.
	
\begin{figure}[htbp!]
		\centering
		\subfloat[]{\includegraphics[width=0.4\textwidth]{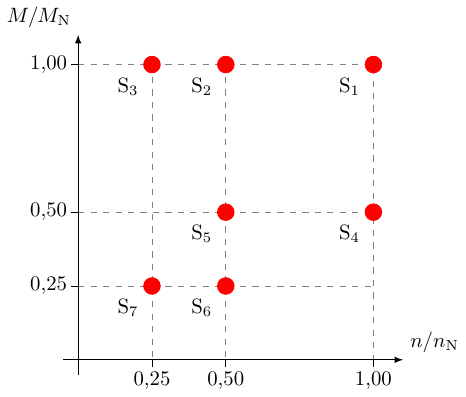}}
		\hfil
		\centering
		\subfloat[]{\includegraphics[width=0.4\textwidth]{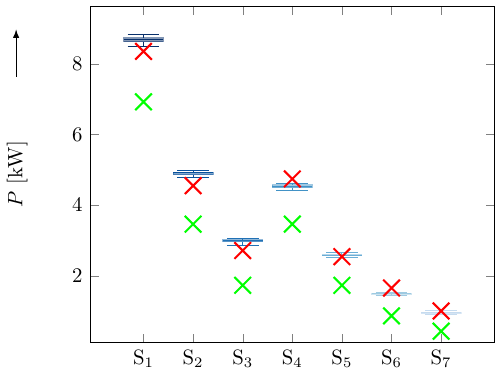}}
		\caption{Panel (a): Standardized reference points for experimental measurements of the grid connection power and power output of the motor. The angular momentum is denoted by $M$ (nominal momentum $M_{\mathrm{N}} = 19~\mathrm{Nm}$) and $n$ stands for the rotational speed  (nominal value $n_{\mathrm{N}} = 3480~\mathrm{min}^{-1}$). Panel (b): Comparison of the experimental values (boxplots) of the grid connection power with the values according to the employed power flow model~\cite{Schutzhold} (red crosses). The green crosses indicate the power output of the motor equal to $2\pi nM$.}
		\label{X4}
	\end{figure}

	In addition, the employed power flow model was experimentally validated at some standardized measurement points: Figure~\ref{X4} presents some results concerning the motor as key part. The simulated grid power values (red crosses in panel b) were compared with the experimentally measured values (boxplots). Additionally, the repective motor power output is shown (green crosses), which allows conclusions about the motor's efficiency ratio (about 81~\% for measurement point~S$_1$).
	
	Table~\ref{Table1} lists boundary values entering (\ref{3.3}~-~\ref{3.5}).
	\section{Numerical Results} \label{Section3}
	This section is divided into two parts: up-travels and down-travels, and both parts contain a systematic evaluation of the optimal trajectories w.r.t.~(\ref{3.1}, \ref{3.11}) over the $x$-$y$ plane covering distance ranges $0.5$~m $\leq s_x \leq 30$~m and $0.5$~m $\leq s_y \leq 20$~m. Therefore, we consider movements with a duration up to 24~s and we are faced with average powers of 30~kW. The results for up-travels can be summarized as follows: 
	\begin{figure}[!h]
		\centering
		\includegraphics[width=0.6\textwidth]{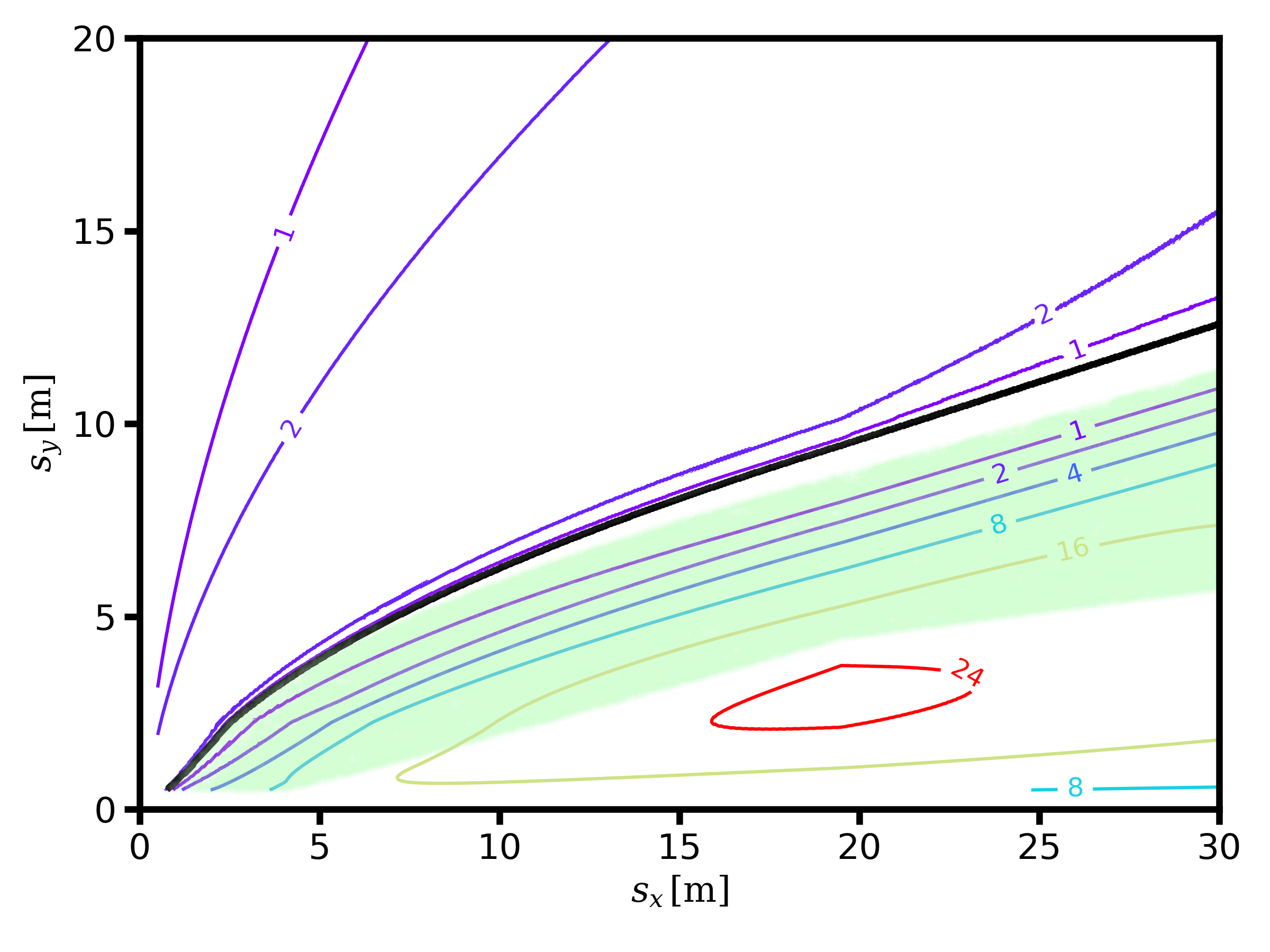} 
		\caption{Contour plot of the energy saving in~\% comparing the energy optimal movement with the fully time minimal movement for up travels. The black curve separates the cases, where the movement is dominated by the running gear (below) or the lifting gear (above). The light green area indicates those combinations of distances $s_x$ and $s_y$, where the energy optimal trajectory can be approximated fairly well by a velocity profile achieving the minimal constant velocity. The island of optimality surrounded by the red curve is formed because of two effects: on the one hand, a longer horizontal movement allows more flexibility for a better recuperation, but a longer travel needs more energy on the other hand.}\label{Figure_7_p}
	\end{figure}
	\begin{itemize}	
		\item The optima of (\ref{3.1}) and (\ref{3.11}) coincide, i.e. the trajectory with optimal usage of recuperation is exactly that one with the least energy consumption.
		\item All optimal trajectories consists of at most seven time intervals, i.e. a time grid $(t_i)_{i=0}^7$ is sufficient. This fact reveals a significant advantage of indirect methods against direct methods which need time step sizes of several milliseconds and therefore grids of hundreds of time points to adequately describe the trajectory.
		\item In addition, there are only two kinds of trajectories: first, the trajectory has only intervals with (distinct) CD type dynamics. Second, the middle interval has EL type dynamics which is surrounded by intervals with CD dynamics (see Figure~\ref{Figure_7_p} below).
	\end{itemize}
	\begin{table}
		\caption{Kinematic parameters $v_{\max}$, $a_{\max}$, $j_{\max}$ defining the velocity, acceleration and jerk bounds, respectively, and entering the optimization problem in (\ref{3.3}~-~\ref{3.5}).}\label{Table1}
		\begin{center}
			\begin{tabular}{lcc}
				\hline
				parameter &\multicolumn{1}{c}{running gear} &\multicolumn{1}{c}{lifting gear} \\ \hline
				$v_{\max}$ &3.0~m/s &0.9~m/s\\
				$a_{\max}$ &0.5~m/s$^2$ &0.6~m/s$^2$\\
				$j_{\max}$ &1.0~m/s$^3$ &0.6~m/s$^3$\\ \hline
			\end{tabular}	
		\end{center}
	\end{table}
	Figure~\ref{Figure_7_p} displays the rate of energy saving, where the following two cases are compared:
	\begin{enumerate}[label=\Alph*)]
		\item The value $E_{\mathrm{rec}}(v)$ if $v$ is chosen as optimal trajectory.
		\item The value $E_{\mathrm{rec}}(v)$ if $v$ is the trajectory of the time-minimal movement.
	\end{enumerate}
	That means that the displayed energy saving rate measures the amount of saved energy (relative to case B) if the trajectory of case A is chosen instead of the trajectory of case B.
	
	The contour plot of Figure~\ref{Figure_7_p} (as well as Figures~\ref{FigureX1} and~\ref{7_n2}) are based on 540~000 trajectories computed. Overall the average value of $E_{\mathrm{rec}}$ is 455.7~kJ and the average saving rate is 5.65~\%. 
	
	The situation changes for downward movements. The main reason for this is the excess energy due to the potential energy of the lifting gear and the fact that the objective function~(\ref{3.1}) penalizes unused energy, whereas objective function~(\ref{3.11}) rewards it. Analogous to Figure~\ref{Figure_7_p}, Figures~\ref{FigureX1} and~\ref{7_n2} show the results for the objective functions energy recuperation~(\ref{3.1}) and energy consumption~(\ref{3.11}), respectively. Let us first turn to the energy recuperation:
	\begin{figure}
		\centering
		\includegraphics[width=0.65\textwidth]{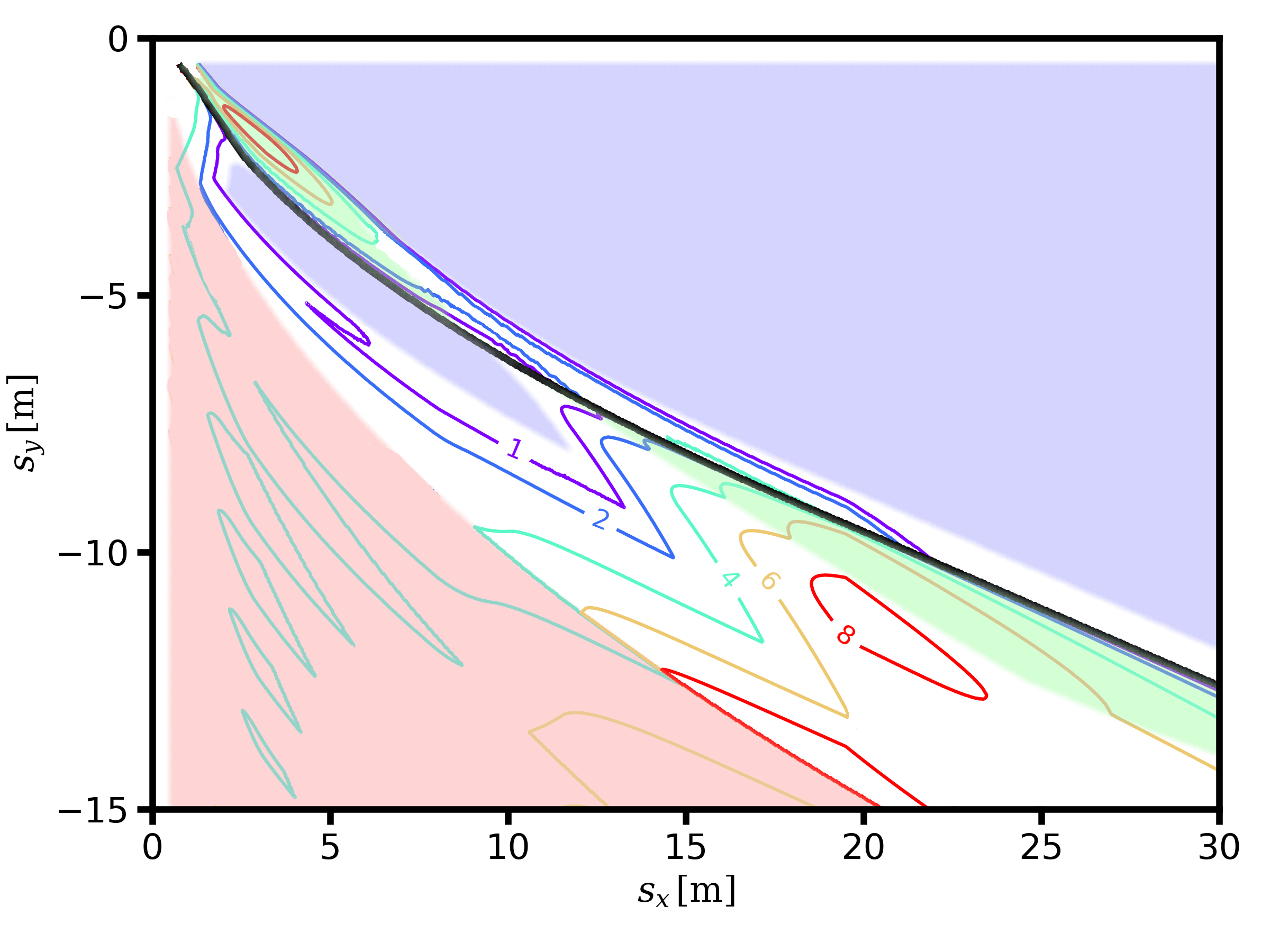}
		\caption{Contour plot of the energy saving in \% if the energy recuperation~(\ref{3.1}) is optimized. In addition to the coloring of Figure~\ref{Figure_7_p}, the blue area indicates those cases, where the energy optimal trajectory is given by a time-minimal movement of both drives. The red area indicates movement where the running gear would start at least twice.\label{FigureX1}}
	\end{figure}
		\begin{figure}
		\centering
		\includegraphics[width=0.7\textwidth]{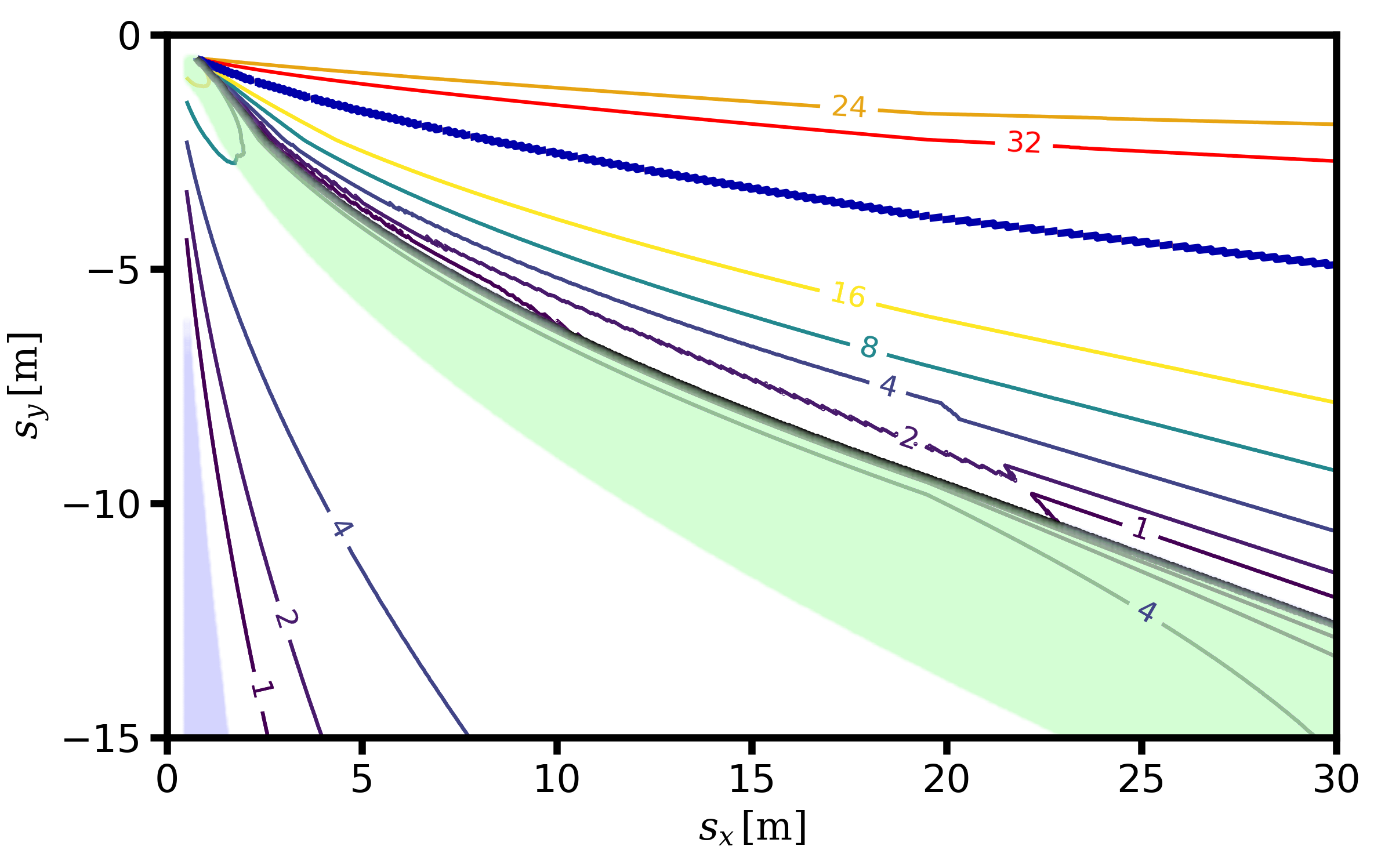}
		\caption{Analogous contour plot (see Fig.~\ref{FigureX1}) for energy consumption~(\ref{3.11}) instead of the recuperation. The thick blue curve separates cases, where $E_{\mathrm{con}}>0$ (above) or $E_{\mathrm{con}}<0$ (below), i.e. all travels located below this curve generate in total an energy excess.\label{7_n2}}
	\end{figure}
	\begin{itemize}
		\item If the movement is dominated by the running gear, the optimal solution corresponds to a time-minimal travel of the lifting gear, because the released potential energy is used for the energy-intensive start-up of the running gear (blue area in Fig.~\ref{FigureX1}).
		\item For long down-travels (dominated by the lifting gear), the large excess energy must be consumed by the running gear which is best achieved in accordance with the objective function~(\ref{3.1}) by having the running gear start and stop several times. The red area in Figure~\ref{FigureX1} indicates those cases in which the running gear starts at least twice. Such a behavior is not only highly detrimental due to increased wear and tear and the potential vibrations of the stacker crane. As a further undesirable side effect, it also generates larger fluctuations in grid power.
		\item From the numerical perspective, these trajectories are the ones with the most time intervals. Time grids with up to 40 grid points are required, which is more than for upward movements, but still significantly fewer than with direct methods.
	\end{itemize}
	This shows that optimizing the energy recuperation alone does not always seem technically reasonable. Another perspective arises if a larger amount of regenerated energy is  unimportant or even advantageous (objective function~(\ref{3.11}), see Figure~\ref{7_n2}). To minimize~(\ref{3.11}), the running gear drives movements with as low accelerations and velocities as possible (for example, the cases of the green area have no acceleration for as long as possible, while the cases in the blue area have no movement at all for as long as possible). In the sense of the objective function~(\ref{3.11}), the time-minimal movement lies not as far away from the optimum for long down-travels as for up-travels, so that the energy saving rates for movements dominated by the lifting gear are relatively small compared to the time-minimal movement.
	
	Once such ``energy maps'' as Figure~\ref{Figure_7_p},~\ref{FigureX1} and~\ref{7_n2} are computed, they can be utilized to develop strategies for stacker crane operations.

	\section{Summary and Outlook} \label{Section4}
	The present case study combines insights gained from a variational principle with numerical tools for continuous convex optimization. Such approaches are refered as indirect methods. The targeted implementation enhances the computations significantly since problem characteristics enter the numerics. Once the analytics and the numerics are set up, it is possible to compute several hundreds of trajectories per second which both ensures the real-time capability and allows to create sufficiently large sets of optimal trajectories serving as data basis for further optimization tasks, e.g. disposition planning. In addition, methods based on a priori discretizations suffer from instabilities caused by the fact that such a discretization does not know anything about important grid points.
	
	In addition to the methodological aspects, the calculations performed showed that maximizing the recuperation, especially for down-travels, leads to technically undesirable solutions. In summary, it can be said that the sole objective of maximizing the recuperation is not practically viable. Another point worth mentioning is that the optimal trajectories with regard to recuperation and energy consumption are similar for up-travels but different for down-travels. From a practical perspective, the finding that trajectories with a minimal constant velocity often prove to be an energetically favourable solution seems relevant as well.
	
	An important step toward a significantly improved energy management for stacker crane opertations in warehouses would be the direct and adaptive connection of trajectory optimization and disposition planning -- a challenge left for future work.
	\section*{Acknowledgement}
	The work of Rico Zöllner was supported by DFG project no.~430149671. The authors thank Frank Schulze for his valuable comments.

\end{document}